\documentclass[a4paper,11pt]{confer}
\usepackage{amsmath,amsfonts,times}
\newcommand{\thefpage}{1}
\newcommand{\thelpage}{1}

\begin{document}

\makeatletter
\@addtoreset{figure}{section}
\def\thefigure{\thesection.\@arabic\c@figure}
\def\fps@figure{h, t}
\@addtoreset{table}{bsection}
\def\thetable{\thesection.\@arabic\c@table}
\def\fps@table{h, t}
\@addtoreset{equation}{section}
\def\theequation{\thesubsection.\arabic{equation}}
\makeatother
\newtheorem{thm}{Theorem}[section]
\newtheorem{prop}[thm]{Proposition}
\newtheorem{lema}[thm]{Lemma}
\newtheorem{cor}[thm]{Corollary}
\newtheorem{defi}[thm]{Definition}
\newtheorem{convention}[thm]{Convention}
\newtheorem{rk}[thm]{Remark}
\newtheorem{exempl}{Example}[section]
\newenvironment{exemplu}{\begin{exempl}  \em}{\hfill $\surd$
\end{exempl}}
\newcommand{\comment}[1]{\par\noindent{\raggedright\texttt{#1}
\par\marginpar{\textsc{Comment}}}}
\newcommand{\todo}[1]{\vspace{5 mm}\par \noindent \marginpar{\textsc{ToDo}}\framebox{\begin{minipage}[c]{0.95 \textwidth}
\tt #1 \end{minipage}}\vspace{5 mm}\par}
\newcommand{\ea}{\mbox{{\bf a}}}
\newcommand{\eu}{\mbox{{\bf u}}}
\newcommand{\ueu}{\underline{\eu}}
\newcommand{\ueo}{\overline{u}}
\newcommand{\oeu}{\overline{\eu}}
\newcommand{\ew}{\mbox{{\bf w}}}
\newcommand{\ef}{\mbox{{\bf f}}}
\newcommand{\eF}{\mbox{{\bf F}}}
\newcommand{\eC}{\mbox{{\bf C}}}
\newcommand{\en}{\mbox{{\bf n}}}
\newcommand{\eT}{\mbox{{\bf T}}}
\newcommand{\eL}{\mbox{{\bf L}}}
\newcommand{\eR}{\mbox{{\bf R}}}
\newcommand{\eV}{\mbox{{\bf V}}}
\newcommand{\eU}{\mbox{{\bf U}}}
\newcommand{\ev}{\mbox{{\bf v}}}
\newcommand{\eve}{\mbox{{\bf e}}}
\newcommand{\uev}{\underline{\ev}}
\newcommand{\eY}{\mbox{{\bf Y}}}
\newcommand{\eK}{\mbox{{\bf K}}}
\newcommand{\eP}{\mbox{{\bf P}}}
\newcommand{\eS}{\mbox{{\bf S}}}
\newcommand{\eJ}{\mbox{{\bf J}}}
\newcommand{\eB}{\mbox{{\bf B}}}
\newcommand{\eH}{\mbox{{\bf H}}}
\newcommand{\leb}{\mathcal{ L}^{n}}
\newcommand{\eI}{\mathcal{ I}}
\newcommand{\eE}{\mathcal{ E}}
\newcommand{\hen}{\mathcal{H}^{n-1}}
\newcommand{\eBV}{\mbox{{\bf BV}}}
\newcommand{\eA}{\mbox{{\bf A}}}
\newcommand{\eSBV}{\mbox{{\bf SBV}}}
\newcommand{\eBD}{\mbox{{\bf BD}}}
\newcommand{\eSBD}{\mbox{{\bf SBD}}}
\newcommand{\ecs}{\mbox{{\bf X}}}
\newcommand{\eg}{\mbox{{\bf g}}}
\newcommand{\paromega}{\partial \Omega}
\newcommand{\gau}{\Gamma_{u}}
\newcommand{\gaf}{\Gamma_{f}}
\newcommand{\sig}{{\bf \sigma}}
\newcommand{\gac}{\Gamma_{\mbox{{\bf c}}}}
\newcommand{\deu}{\dot{\eu}}
\newcommand{\dueu}{\underline{\deu}}
\newcommand{\dev}{\dot{\ev}}
\newcommand{\duev}{\underline{\dev}}
\newcommand{\weak}{\stackrel{w}{\approx}}
\newcommand{\mild}{\stackrel{m}{\approx}}
\newcommand{\strong}{\stackrel{s}{\approx}}
\newcommand{\weakdown}{\rightharpoondown}
\newcommand{\opg}{\stackrel{\mathfrak{g}}{\cdot}}
\newcommand{\opunu}{\stackrel{1}{\cdot}}
\newcommand{\opdoi}{\stackrel{2}{\cdot}}
\newcommand{\opn}{\stackrel{\mathfrak{n}}{\cdot}}
\newcommand{\opx}{\stackrel{x}{\cdot}}
\newcommand{\tr}{\ \mbox{tr}}
\newcommand{\Ad}{\ \mbox{Ad}}
\newcommand{\ad}{\ \mbox{ad}}

\title*{Dilatation structures in sub-riemannian geometry}
\author{MARIUS BULIGA \email{Marius.Buliga@imar.ro}}
\institute{Institute of Mathematics, Romanian Academy \\
P.O. BOX 1-764, RO 014700\\
Bucure\c sti, Romania}

\keywords{Carnot groups, dilatation structures, sub-riemannian manifolds}
\amssubj{53C17, 22E20, 20F65}
\maketitle

\begin{abstract}
Based on the notion of dilatation structure \cite{buligadil1}, we give
an intrinsic treatment to sub-riemannian geometry, started
 in the paper \cite{buligaradon}.
  Here we prove that  regular sub-riemannian manifolds admit dilatation
structures. From the existence of normal frames
proved by Bella\"{\i}che we deduce the rest of the properties of  regular
sub-riemannian manifolds by using the formalism of dilatation structures.
\end{abstract}

\section{Introduction}

Sub-riemannian geometry is the modern incarnation of non-holonomic spaces,
discovered in 1926 by the romanian mathematician Gheorghe Vr\u anceanu \cite{vra1},
\cite{vra2}. The sub-riemannian geometry is the study of non-holonomic spaces
endowed with a Carnot-Carath\'eodory distance. Such spaces appear in applications
to thermodynamics, to the mechanics of non-holonomic systems, in the study of
hypo-elliptic operators cf. H\"ormander \cite{hormander}, in harmonic analysis on
homogeneous cones cf. Folland, Stein \cite{folstein}, and
as boundaries of CR-manifolds.

The interest in these spaces comes from several intriguing features which they have:
from the metric point of view they are fractals (the Hausdorff dimension with
respect to the Carnot-Carath\'eodory distance is strictly bigger than the
topological dimension, cf. Mitchell \cite{mit}); the metric tangent space to a point of a regular
 sub-riemannian manifold is a Carnot
group (a simply connected nilpotent Lie group with a positive graduation), also
known classicaly as a homogeneous cone; the asymptotic space (in the sense of
Gromov-Hausdorff distance) of a finitely generated group  with polynomial
growth is also a Carnot group, by a famous theorem of Gromov \cite{gromovgr} wich
leads to an inverse to the Tits alternative; finally, on such spaces we have
enough structure to develop a differential calculus resembling to the one
proposed by  Cheeger \cite{cheeger} and to prove theorems like Pansu' version
of Rademacher theorem \cite{pansu}, leading to an ingenious proof of a Margulis
rigidity result.

There are several fundamental papers dedicated to the establisment of the
sub-riemannian geometry, among them Mitchell \cite{mit},
Bella\"{\i}che \cite{bell}, a substantial paper of Gromov asking for an
 intrinsic point of view for sub-riemannian geometry  \cite{gromovsr}, Margulis,
 Mostow \cite{marmos1}, \cite{marmos2}, dedicated to Rademacher theorem for
 sub-riemannian manifolds and to the construction of a tangent bundle of such
 manifolds, and Vodopyanov \cite{vodopis} (among other papers), concerning the same
 subject.

 There is a reason for the existence of so many papers, written by important
 mathematicians, on the same subject: the fundamental geometric
 properties of sub-riemannian manifolds are very difficult to prove. Maybe
 the most difficult problem is to provide a rigorous construction of the
 tangent bundle of such a manifold, starting from the properties of the
 Carnot-Carath\'eodory distance, and somehow allowing to generalize Pansu'
 differential calculus.

 In several articles devoted to sub-riemannian geometry,
these fundamental results  were proved using
differential geometry tools, which are  not intrinsic to
sub-riemannian geometry, therefore leading to very long proofs, sometimes with
unclear parts, corrected or clarified in other papers dedicated to the same subject. 

The fertile ideas of Gromov, Bella\"{\i}che and other founders of the field of analysis 
in sub-riemannian spaces are now developed into a hot research area. For the study of 
sub-riemannian geometry under weaker than usual regularity hypothesis see for example 
the string of papers by Vodopyanov, among them \cite{vodopis}, \cite{vodopis2}. In these 
papers  Vodopyanov constructs a tangent bundle structure for a sub-riemannian 
manifold, under weak regularity hypothesis, by using notions as horizontal convergence. 

Based on the notion of dilatation structure \cite{buligadil1}, I tried to give a
an intrinsic treatment to sub-riemannian geometry in the paper \cite{buligaradon},
after a series of articles \cite{buliga2}, \cite{buliga3}, \cite{buliga4} dedicated
to the sub-riemannian geometry of Lie groups endowed with left invariant
distributions.

In this article we show that normal frames are the central objects in the establishment of 
fundamental properties in sub-riemannian geometry, in the following precise sense. 
 We prove  that for  regular sub-riemannian manifolds, 
the existence of normal frames 
(definition \ref{defnormal}) implies that induced   dilatation
structures exist (theorems \ref{thma3}, \ref{pcomp3}).  The existence of normal frames has been proved 
by Bella\"{\i}che  \cite{bell}, starting
with theorem 4.15 and ending in the first half of section 7.3 (page 62). From these facts all 
 classical results concerning
the structure of the tangent space to a point of a regular sub-riemannian manifold
can be deduced as straightforward consequences of the structure theorems \ref{thcone1},
\ref{thcone}, \ref{tgene}, \ref{cortang} from 
the formalism of dilatation structures. 

In conclusion, our purpose is twofold: (a) we try to show that basic results in sub-riemannian geometry 
are particular cases  of the abstract theory of dilatation structures, and (b) we try to minimize 
the contribution of classical differential calculus in the proof of these basic results, by showing 
that in fact the differential calculus on the sub-riemannian manifold is needed only for proving 
that normal frames exist and after this stage an intrinsic way of reasoning is possible. 

If we take the point of view of Gromov, that the only intrinsic object on a sub-riemannian  
manifold is the Carnot-Carath\'eodory  distance, the underlying differential structure of the 
manifold is clearly not intrinsic. Nevertheless in all proofs that I know this 
differential structure is heavily used. Here we  try to prove that in fact it is  
sufficient to  take as intrinsic objects of sub-riemannian geometry 
the Carnot-Carath\'eodory distance and dilatation structures compatible with it. 

The closest results along these lines are maybe the ones of Vodopyanov. There is a clear correspondence 
between his way of defining the tangent bundle of a sub-riemannian manifold and the way of dilatation 
structures. In both cases the tangent space to a point is defined only locally, as a neighbourhood 
of the point, in the manifold, endowed with a local group operation. Vodopyanov proves the existence 
of the (locally defined) operation under very weak regularity assumptions on the sub-riemannian 
manifold. The main tool of his proofs is nevertheless the differential structure of the underlying manifold. 
In distinction, we prove in \cite{buligadil1}, in  an abstract setting, 
  that the very existence of a dilatation structure induces a locally defined operation. Here we show 
that the differential structure of the underlying manifold  is important only in order to prove 
 that dilatation structures can indeed be constructed  from normal frames.

{\bf Aknowledgements.} I want to thank the anonymous referee for several 
valuable suggestions, as mentioned further in the paper. 

\section{Metric profiles}

\paragraph{Notations.} The space $CMS$ is  the collection of isometry classes of pointed compact metric
spaces. The notation used for elements of $CMS$ is of the type  $[X,d,x]$,
representing the
equivalence class of the pointed compact metric space $(X,d,x)$ with respect
to (pointed) isometry. The open ball of radius $r>0$ and center $x \in (X,d)$ is
denoted by $B(x,r)$ or $\displaystyle B_{d}(x,r)$ if we want to emphasize the
dependence on the distance $d$. The notation for a closed ball is obtained by
adding an overline to the notation for the open ball.   The distance on $CMS$ is the Gromov-Hausdorff distance $\displaystyle
d_{GH}$ between
(isometry classes of) pointed metric spaces and the topology is induced by
this distance. For the Gromov-Hausdorff distance see Gromov \cite{gromov}. We denote by $\mathcal{O}(\varepsilon)$  a positive function such that
$\displaystyle \lim_{\varepsilon \rightarrow 0} \mathcal{O}(\varepsilon) \ = \ 0$.

To any locally compact metric  space there is an associated   metric profile
(Buliga \cite{buliga3}, \cite{buliga4}).

\begin{defi}
 The metric profile associated to the
 locally compact metric space $(M,d)$ is  the assignment (for small enough $\varepsilon > 0$)
$$(\varepsilon > 0 , \ x \in M) \ \mapsto \  \mathbb{P}^{m}(\varepsilon, x) = [\bar{B}(x,1) , \frac{1}{\varepsilon} d, x] \in CMS$$
\label{dmprof}
\end{defi}

We may define a notion of metric profile which is more general than the previous
one.

\begin{defi}
\label{dprofile}
A metric profile is a curve $\mathbb{P}:[0,a] \rightarrow CMS$ such that
\begin{enumerate}
\item[(a)] it is continuous at $0$,
\item[(b)]for any $\mu \in [0,a]$ and  $\varepsilon \in (0,1]$ we have
$$d_{GH} (\mathbb{P}(\varepsilon \mu),
\mathbb{P}^{m}_{d_{\mu}}(\varepsilon,x_{\mu}))
\  = \ O(\mu)$$
\end{enumerate}
The function $\mathcal{O}(\mu)$ may change with  $\varepsilon$.
We used the notations
$$
\mathbb{P}(\mu) = [\bar{B}(x,1) ,d_{\mu}, x_{\mu}] \quad \mbox{  and } \quad
\mathbb{P}^{m}_{d_{\mu}}(\varepsilon,x) =
\left[\bar{B}(x,1),\frac{1}{\varepsilon}d_{\mu}, x_{\mu}\right]
$$
\end{defi}

We shall unfold further the definition \ref{dprofile} in order to clearly
understand what is a metric profile.  For any $\mu \in (0,a]$ and for any
$b > 0$ there is $\varepsilon(\mu, b) \in (0,1)$ such that for any
$\varepsilon \in (0,\varepsilon(\mu,b))$ there exists a relation
$\displaystyle \rho = \rho_{\varepsilon, \mu} \subset \bar{B}_{d_{\mu}}(x_{\mu},
\varepsilon) \times \bar{B}_{d_{\mu \varepsilon}}(x_{\mu \varepsilon}, 1)$ such
that:
\begin{enumerate}
\item[1.] $\displaystyle dom \ \rho_{\varepsilon, \mu}$ is $b$-dense in
$\displaystyle \bar{B}_{d_{\mu}}(x_{\mu}, \varepsilon)$,
\item[2.] $\displaystyle im \ \rho_{\varepsilon, \mu}$ is $b$-dense in
$\displaystyle \bar{B}_{d_{\mu \varepsilon}}(x_{\mu \varepsilon}, 1)$,
\item[3.] $\displaystyle (x_{\mu}, x_{\mu \varepsilon}) \in \rho_{\varepsilon,
\mu}$,
\item[4.] for all $\displaystyle x,y \in \ dom \ \rho_{\varepsilon, \mu}$ we have
 $\displaystyle \mid \frac{1}{\varepsilon} d_{\mu}(x,y) - d_{\mu
 \varepsilon}( x', y')   \mid \
 \leq \  b$,
 for any $x', y'$ such that $\displaystyle (x,x'), (y,y') \in \rho_{\varepsilon,
 \mu}$.
\end{enumerate}

Therefore a metric profile gives two types of information:
\begin{enumerate}
\item[-] a distance estimate like the one  from point 4 above,
\item[-] an "approximate shape" estimate, like in the points 1--3, where we
see that two sets, namely the balls
$\displaystyle \bar{B}_{d_{\mu}}(x_{\mu}, \varepsilon)$ and
$\displaystyle \bar{B}_{d_{\mu \varepsilon}}(x_{\mu \varepsilon}, 1)$, are
approximately isometric.
\end{enumerate}

The metric profile $\varepsilon \mapsto \mathbb{P}^{m}(\varepsilon, x)$ of a metric space $(M,d)$ for  a fixed $x \in M$ is a metric profile in the sense of the definition \ref{dprofile} if and only if the space
$(M,d)$ admits a tangent space in $x$. Here is the general definition of a
tangent space in the metric sense.

\begin{defi}
\label{defmetspace}
A (locally compact) metric space $(M,d)$ admits a (metric) tangent space in $x \in M$ if the associated metric profile
$\varepsilon \mapsto \mathbb{P}^{m}(\varepsilon, x)$ (as in definition \ref{dmprof})  admits a prolongation by continuity in
$\varepsilon = 0$, i.e if the following limit exists:
\begin{equation}
\label{limmetspace}
[T_{x}M,d^{x}, x]  =  \lim_{\varepsilon \rightarrow 0} \mathbb{P}^{m}(\varepsilon, x)
\end{equation}
\end{defi}

Metric tangent spaces are   metric cones.

\begin{defi}
\label{defmetcone}
A metric cone $(X,d, x)$ is a locally compact metric space $(X,d)$, with a marked point $x \in X$ such
that for any $a,b \in (0,1]$ we have
\[\displaystyle \mathbb{P}^{m}(a,x)  =  \mathbb{P}^{m}(b,x)\]
\end{defi}

Metric cones have the  simplest metric profile, which  is one with the property:
$\displaystyle (\bar{B}(x_{b}, 1), d_{b}, x_{b}) = (X, d_{b}, x)$.
In particular metric cones have dilatations.

\begin{defi}
Let $(X,d, x)$ be a metric cone. For any $\varepsilon \in (0,1]$  a dilatation is a function $\displaystyle \delta^{x}_{\varepsilon}: \bar{B}(x,1) \rightarrow \bar{B}(x,\varepsilon)$ such that
\begin{enumerate}
\item[(a)]
$\displaystyle \delta^{x}_{\varepsilon}(x) = x$,
\item[(b)]
for any $u,v \in X$ we have
\[d\left(\delta^{x}_{\varepsilon}(u), \delta^{x}_{\varepsilon}(v)\right) =\varepsilon \, d(u,v) \]
\end{enumerate}
\end{defi}

The existence of dilatations for metric cones comes from the definition
\ref{defmetcone}. Indeed, dilatations are just  isometries from
$\displaystyle (\bar{B}(x,1), d, x)$ to $ (\bar{B}, \frac{1}{a}d, x)$.

\section{Sub-riemannian manifolds}

Let $M$ be a connected $n$ dimensional real manifold. A distribution
is a smooth  subbundle $D$ of $M$. To any point
$x \in M$ there is associated the vector space $D_{x} \subset T_{x}M$.
The dimension of the distribution $D$ at point $x \in M$ is denoted by
$$m(x) = \, dim \, D_{x}$$
The distribution is smooth, therefore the function $x \in M \mapsto m(x)$ is
locally constant. We suppose further that the dimension of the distribution is
globally constant and we  denote it by $m$ (thus $m = m(x)$ for any
$x \in M$). Clearly $m \leq n$; we are interested in the case $m < n$.

A horizontal curve $c:[a,b] \rightarrow M$ is a curve which is almost everywhere derivable and for
almost any $t \in [a,b]$ we have $\displaystyle \dot{c}(t) \in D_{c(t)}$.
The class of horizontal curves will be denoted by $Hor(M,D)$.

Further we shall use the following notion of non-integrability of the
distribution $D$.

\begin{defi}
The distribution $D$ is completely non-integrable if $M$ is locally connected
by horizontal curves curves  $c \in Hor(M,D)$.
\label{defcnint}
\end{defi}

A sufficient condition for the distribution $D$ to be completely non-integrable
is given by Chow condition (C) \cite{chow}.

\begin{thm} (Chow) Let $D$ be a distribution of dimension $m$  in the manifold
$M$. Suppose there is a positive integer number $k$ (called the rank of the
distribution $D$) such that for any $x \in X$ there is a topological  open ball
$U(x) \subset M$ with $x \in U(x)$ such that there are smooth vector fields
$\displaystyle X_{1}, ..., X_{m}$ in $U(x)$ with the property:

(C) the vector fields $\displaystyle X_{1}, ..., X_{m}$ span $\displaystyle
D_{x}$ and these vector fields together with  their iterated
brackets of order at most $k$ span the tangent space $\displaystyle T_{y}M$
at every point $y \in U(x)$.

Then the distribution $D$ is completely non-integrable in the sense of
definition \ref{defcnint}.
\label{tchow}
\end{thm}

\begin{defi}
A sub-riemannian (SR) manifold is a triple $(M,D, g)$, where $M$ is a
connected manifold, $D$ is a completely non-integrable distribution on $M$, and $g$ is a metric (Euclidean inner-product) on the distribution (or horizontal bundle)  $D$.
\label{defsr}
\end{defi}

\subsection{The Carnot-Carath\'eodory distance}

Given a distribution $D$ which satisfies the hypothesis of Chow theorem
\ref{tchow}, let us consider a point   $x \in M$,  its neighbourhood $U(x)$,
and the vector fields $\displaystyle X_{1}, ..., X_{m}$ satisfying the
condition (C).

One can define on $U(x)$ a filtration of bundles
as follows. Define first the class of horizontal vector fields on $U$:
$$\mathcal{X}^{1}(U(x),D) \ = \ \left\{ X \in \mathcal{X}^{\infty}(U) \mbox{ : }
\forall y \in U(x) \ , \ X(y) \in D_{y} \right\}$$
Next, define inductively for all positive integers $j$:
$$ \mathcal{X}^{j+1} (U(x),D) \ = \ \mathcal{X}^{j}(U(x),D) \, + \,
 [ \mathcal{X}^{1}(U(x),D),
\mathcal{X}^{j}(U(x),D)]$$
Here $[ \cdot , \cdot ]$ denotes the bracket of vector fields.
We obtain therefore a filtration $\displaystyle \mathcal{X}^{j}(U(x),D) \subset
\mathcal{X}^{j+1} (U(x),D)$.
Evaluate now this filtration at $y \in U(x)$:
$$V^{j}(y,U(x),D) \ = \ \left\{ X(y) \mbox{ : } X \in \mathcal{X}^{j}(U(x),D)\right\}$$
According to Chow theorem there is a positive integer $k$ such that
 for all $y \in U(x)$ we have
$$D_{y} =  V^{1}(y,U(x),D) \subset V^{2}(y,U(x),D) \subset ... \subset
V^{k}(y, U(x),D) = T_{y}M $$
Consequently, to the sub-riemannian manifold is associated the string of
numbers:
$$ \nu_{1}(y) = \dim V^{1}(y, U(x),D) < \nu_{2}(y) = \dim V^{2}(y, U(x),D)
< ... < n = \dim M$$
Generally $k$, $\displaystyle \nu_{j}(y)$
may vary from a point to another.

The number $k$ is called the step of the distribution at
$y $.

\begin{defi}
The distribution $D$ is regular if  $\displaystyle \nu_{j}(y)$ are constant on
the manifold $M$.
The sub-riemannian manifold $M,D,g)$ is regular if $D$ is regular and for any
$x \in M$ there is a topological ball $U(x) \subset M$ with $x \in U(M)$ and
 an orthonormal (with respect to the metric $g$)  family of
  smooth vector fields $\displaystyle \left\{ X_{1}, ..., X_{m}\right\}$ in
  $U(x)$ which satisfy the condition (C).
\label{dreg}
\end{defi}

The lenght of  a horizontal curve is
$$l(c) \ = \ \int_{a}^{b} \left(g_{c(t)} (\dot{c}(t),
\dot{c}(t))\right)^{\frac{1}{2}} \mbox{ d}t$$
The length depends on the metric $g$.

\begin{defi}
The Carnot-Carath\'eodory distance (or CC distance) associated to the sub-riemannian manifold is the
distance induced by the length $l$ of horizontal curves:
$$d(x,y) \ = \ \inf \left\{ l(c) \mbox{ : } c \in Hor(M,D) \
, \ c(a) = x \ , \  c(b) = y \right\} $$
\end{defi}

The Chow theorem ensures the existence of a horizontal path linking any two sufficiently
closed points, therefore the CC distance is  locally finite. The distance
depends only on the distribution $D$ and metric $g$, and not on the choice of
vector fields $\displaystyle X_{1}, ..., X_{m}$ satisfying the condition (C).
The space $(M,d)$ is locally compact and complete, and the topology
induced by the distance $d$ is the same as the topology of the manifold $M$.
(These important details may be recovered from reading carefully the
constructive proofs of Chow theorem given by  Bella\"{\i}che
\cite{bell} or Gromov \cite{gromovsr}.)

\subsection{Normal frames}

In the following we stay in a small open neighbourhood of an arbitrary, but
fixed point $\displaystyle x_{0} \in M$. All results are local in
nature (that is they hold for some small open neighbourhood of an arbitrary, but
fixed point of the manifold $M$). That is why we shall no longer mention the
dependence of various objects on $\displaystyle x_{0}$, on the neighbourhood
$\displaystyle U(x_{0})$, or the distribution $D$.

We shall work further only with regular sub-riemannian manifolds, if not
otherwise stated. The topological dimension of $M$ is denoted by $n$, the
 step of the regular sub-riemannian manifold $(M,D,g)$ is denoted by $k$, the
dimension of the distribution is $m$, and there are numbers $\displaystyle
\nu_{j}$, $j = 1, ..., k$ such that for any $x \in M$ we have
$\displaystyle dim \, V^{j}(x) = \nu_{j}$. The Carnot-Carath\'eodory distance
is denoted by $d$.

\begin{defi}
An adapted frame $\displaystyle \left\{ X_{1}, ... , X_{n} \right\}$ is a
collection of smooth vector fields which is obtained by the construction
described below.

We start with a collection $X_{1}, ... , X_{m}$ of vector fields which satisfy
the condition (C). In particular  for any point $x$ the vectors  $\displaystyle
X_{1}(x), ... , X_{m}(x)$ form a basis for $\displaystyle D_{x}$.  We further
associate to any word $\displaystyle a_{1} .... a_{q}$ with letters in the
alphabet $1, ... ,m$ the  multi-bracket
$\displaystyle [X_{a_{1}}, [ ... , X_{a_{q}}] ... ]$.

One can add,  in the lexicographic order, $n-m$ elements to the set
$\displaystyle \left\{ X_{1}, ... , X_{m} \right\}$ until we get a collection
$\displaystyle \left\{ X_{1}, ... , X_{n} \right\}$ such that:
for any $j = 1, ..., k$ and for any point $x$ the set
$\displaystyle \left\{X_{1}(x), ..., X_{\nu_{j}}(x) \right\}$ is a basis for
$\displaystyle V^{j}(x)$.
\label{defadapt}
\end{defi}

Let $\displaystyle \left\{ X_{1}, ... , X_{n} \right\}$ be an adapted frame.
For any $j = 1, ..., n$  the degree $\displaystyle deg \, X_{j}$ of the vector
field $\displaystyle X_{j}$  is defined
as the only positive integer $p$ such that for any point $x$ we have
$$X_{j}(x) \in V^{p}_{x} \setminus V^{p-1}(x)$$

Further we define normal frames. The name has been used by Vodopyanov
\cite{vodopis}, but for a slightly different object. The existence of normal
frames in the sense of the following definition is the hardest technical
problem in the classical establishment of sub-riemannian geometry. For the 
informed reader the referee 
pointed out that condition (a) Definition \ref{defnormal} 
is a part of the conclusion of Gromov approximation theorem, namely when one point coincides 
with the center of nilpotentization; also condition (b) is equivalent with a statement of 
Gromov concerning the convergence of rescaled vector fields to their nilpotentization 
(an informed reader must at least follow in all details the papers Bella\"{\i}che 
\cite{bell} and Gromov \cite{gromovsr}, where differential calculus in the classical 
sense is heavily used). Therefore the conditions of Definition \ref{defnormal} concentrate 
that part of the foundations of sub-riemannian geometry which makes use of classical 
differential calculus. 

The key details in the Definition below are uniform convergence assumptions. This is in line 
with Gromov suggestions in the last section of Bella\"{\i}che \cite{bell}. 

\begin{defi}
An adapted frame $\displaystyle \left\{ X_{1}, ... , X_{n} \right\}$ is a normal
frame if the following two conditions are satisfied:
\begin{enumerate}
\item[(a)] we have the limit
$$\lim_{\varepsilon \rightarrow 0_{+}} \frac{1}{\varepsilon} \, d\left(
\exp \left( \sum_{1}^{n}\varepsilon^{deg\, X_{i}} a_{i} X_{i} \right)(y), y \right) \ = \ A(y, a) \in
(0,+\infty)$$
uniformly with respect to $y$ in compact sets and $\displaystyle a=(a_{1}, ...,
a_{n}) \in W$, with $\displaystyle W \subset \mathbb{R}^{n}$ compact
neighbourhood of $\displaystyle 0 \in \mathbb{R}^{n}$,
\item[(b)] for any compact set $K\subset M$ with diameter (with respect to the
distance $d$) sufficiently small,  and for any $i = 1, ..., n$ there
are functions 
$$ P_{i}(\cdot, \cdot, \cdot): U_{K} \times U_{K} \times K \rightarrow  \mathbb{R}$$
 with $\displaystyle U_{K} \subset \mathbb{R}^{n}$ a sufficiently small compact neighbourhood of 
$\displaystyle 0 \in \mathbb{R}^{n}$ such that for any   $x \in K$ and any $\displaystyle 
a,b \in U_{K}$ we have
$$\exp \left( \sum_{1}^{n} a_{i} X_{i} \right) (x) \ = \
\exp \left( \sum_{1}^{n} P_{i}(a, b, y) X_{i} \right) \circ \exp \left( \sum_{1}^{n}  b_{i} X_{i} \right) (x) $$
and such that the following limit exists
$$\lim_{\varepsilon \rightarrow 0_{+}}
\varepsilon^{-deg \, X_{i}} P_{i}(\varepsilon^{deg \, X_{j}} a_{j}, \varepsilon^{deg \, X_{k}} b_{k}, x)   \in
\mathbb{R}$$
and it is uniform with respect to $x  \in K$ and $\displaystyle a, b \in U_{K}$.
\end{enumerate}
\label{defnormal}
\end{defi}

The existence of normal frames is proven in Bella\"{\i}che \cite{bell}, starting
with theorem 4.15 and ending in the first half of section 7.3 (page 62).

In order to understand normal frames let us look to the
 case of a Lie group $G$ endowed with a left invariant distribution.
 The distribution is completely non-integrable if it is generated by the left
 translation of a vector subspace $D$ of the algebra
$\mathfrak{g} = T_{e}G$ which bracket generates the whole algebra
$\mathfrak{g}$. Take $\displaystyle \left\{ X_{1}, ..., X_{m}\right\}$ a collection of $m = \,
dim \, D$ left invariant independent vector fields and define with their help
an adapted frame, as explained in definition \ref{defadapt}. Then the adapted
frame $\displaystyle \left\{ X_{1}, ..., X_{n}\right\}$ is in fact normal.

\section{Dilatation structures}

In this section we review the definition and main properties
of a  dilatation structure, according to \cite{buligadil1},
\cite{buligacont}.

\subsection{The axioms of a dilatation structure}

Further are listed the  axioms of  a dilatation structure $(X,d,\delta)$,
 starting with axiom 0,
which is a preparation for the axioms which follow.

We restrict the generality from \cite{buligadil1} to the case which is related to sub-riemannian
geometry, that is we shall consider only dilatations $\displaystyle
\delta^{x}_{\varepsilon}$ with $\varepsilon \in (0,+\infty)$.

\begin{enumerate}
\item[{\bf A0.}] The dilatations $$ \delta_{\varepsilon}^{x}: U(x)
\rightarrow V_{\varepsilon}(x)$$ are defined for any
$\displaystyle \varepsilon \in (0,1]$. The sets
$\displaystyle U(x), V_{\varepsilon}(x)$ are open neighbourhoods of $x$.
All dilatations are homeomorphisms (invertible, continuous, with
continuous inverse).

We suppose  that there is a number  $1<A$ such that for any $x \in X$ we have
$$\bar{B}_{d}(x,A) \subset U(x)  \ .$$
 We suppose that for all $\varepsilon \in
(0,1)$, we have
$$ B_{d}(x, \varepsilon) \subset \delta_{\varepsilon}^{x} B_{d}(x,A)
\subset V_{\varepsilon}(x) \subset U(x) \ .$$

 There is a number $B \in (1,A)$ such that  for
 any $\varepsilon \in (1,+\infty)$ the associated dilatation
$$\delta^{x}_{\varepsilon} : W_{\varepsilon}(x) \rightarrow B_{d}(x,B) \ , $$
is injective, invertible on the image. We shall suppose that
$\displaystyle  W_{\varepsilon}(x)$ is a open neighbourhood of $x$,
$$V_{\varepsilon^{-1}}(x) \subset W_{\varepsilon}(x) $$
and that for all $\displaystyle \varepsilon \in (0,1)$ and
$\displaystyle u \in U(x)$ we have
$$\delta_{\varepsilon^{-1}}^{x} \ \delta^{x}_{\varepsilon} u \ = \ u \ .$$
\end{enumerate}

We have therefore  the following string of inclusions, for any
$\varepsilon \in (0,1)$, and any $x \in X$:
$$ B_{d}(x, \varepsilon) \subset \delta^{x}_{\varepsilon}  B_{d}(x, A)
\subset V_{\varepsilon}(x) \subset
W_{\varepsilon^{-1}}(x) \subset \delta_{\varepsilon}^{x}  B_{d}(x, B) \quad . $$

A further technical condition on the sets  $\displaystyle V_{\varepsilon}(x)$ and $\displaystyle W_{\varepsilon}(x)$  will be given just before the axiom A4. (This condition will be counted as part of
axiom A0.)

\begin{enumerate}
\item[{\bf A1.}]  We  have
$\displaystyle  \delta^{x}_{\varepsilon} x = x $ for any point $x$.
We also have $\displaystyle \delta^{x}_{1} = id$ for any $x \in X$.

Let us define the topological space
$$ dom \, \delta = \left\{ (\varepsilon, x, y) \in (0,+\infty) \times X
\times X \mbox{ : } \quad \mbox{ if } \varepsilon \leq 1 \mbox{ then } y
\in U(x) \,
\, ,
\right.$$
$$\left. \mbox{  else } y \in W_{\varepsilon}(x) \right\} $$
with the topology inherited from the product topology on
$(0,+\infty) \times X \times X$. Consider also
$\displaystyle Cl(dom \, \delta)$,
the closure of $dom \, \delta$ in $\displaystyle [0,+\infty) \times X \times X$
with product topology. The function $\displaystyle \delta : dom \, \delta
\rightarrow  X$ defined by $\displaystyle \delta (\varepsilon,  x, y)  =
\delta^{x}_{\varepsilon} y$ is continuous. Moreover, it can be continuously
extended to $\displaystyle Cl(dom \, \delta)$ and we have
$$\lim_{\varepsilon\rightarrow 0} \delta_{\varepsilon}^{x} y \, = \, x \quad . $$

\item[{\bf A2.}] For any  $x, \in K$, $\displaystyle \varepsilon, \mu \in (0,1)$
 and $\displaystyle u \in
\bar{B}_{d}(x,A)$   we have:
$$ \delta_{\varepsilon}^{x} \delta_{\mu}^{x} u  = \delta_{\varepsilon \mu}^{x} u  \ .$$

\item[{\bf A3.}]  For any $x$ there is a  function $\displaystyle (u,v) \mapsto d^{x}(u,v)$, defined for any $u,v$ in the closed ball (in distance d) $\displaystyle
\bar{B}(x,A)$, such that
$$\lim_{\varepsilon \rightarrow 0} \quad \sup  \left\{  \mid
\frac{1}{\varepsilon} d(\delta^{x}_{\varepsilon} u, \delta^{x}_{\varepsilon} v) \ - \ d^{x}(u,v) \mid \mbox{ :  } u,v \in \bar{B}_{d}(x,A)\right\} \ =  \ 0$$
uniformly with respect to $x$ in compact set.

\end{enumerate}

Remark that  $d^{x}$ may be  a degenerated distance: there might exist
$\displaystyle v,w \in U(x)$ such that $\displaystyle d^{x}(v,w) = 0$.

For  the following axiom to make sense we impose a technical condition on the co-domains $\displaystyle V_{\varepsilon}(x)$: for any compact set $K \subset X$ there are $R=R(K) > 0$ and
$\displaystyle \varepsilon_{0}= \varepsilon(K) \in (0,1)$  such that
for all $\displaystyle u,v \in \bar{B}_{d}(x,R)$ and all
$\displaystyle \varepsilon  \in (0,\varepsilon_{0})$,  we have
$$\delta_{\varepsilon}^{x} v \in W_{\varepsilon^{-1}}( \delta^{x}_{\varepsilon}u) \ .$$

With this assumption the following notation makes sense:
$$\Delta^{x}_{\varepsilon}(u,v) = \delta_{\varepsilon^{-1}}^{\delta^{x}_{\varepsilon} u} \delta^{x}_{\varepsilon} v . $$
The next axiom can now be stated:
\begin{enumerate}
\item[{\bf A4.}] We have the limit
$$\lim_{\varepsilon \rightarrow 0}  \Delta^{x}_{\varepsilon}(u,v) =  \Delta^{x}(u, v)  $$
uniformly with respect to $x, u, v$ in compact set.
\end{enumerate}

\begin{defi}
A triple $(X,d,\delta)$ which satisfies A0, A1, A2, A3, but $\displaystyle d^{x}$ is degenerate for some
$x\in X$, is called degenerate dilatation structure.

If the triple $(X,d,\delta)$ satisfies A0, A1, A2, A3, A4 and
 $\displaystyle d^{x}$ is non-degenerate for any $x\in X$, then we call it  a
 dilatation structure.
 \label{defweakstrong}
\end{defi}

\subsection{Metric profile of a dilatation structure}

Here we describe the metric profile associated to a dilatation structure. This
will be relevant further for understanding the geometry of the metric tangent
spaces of regular sub-riemannian manifolds.

The following result is a reformulation of theorem 6 \cite{buligadil1}.

\begin{thm}
\label{thcone1}
Let $(X,d,\delta)$ be a dilatation structure, $x \in X$ a point in $X$, $\mu >
0$ sufficiently
small,
 and let  $(\delta, \mu, x)$ be  the distance on
$\displaystyle \bar{B}_{d^{x}}(x,1) = \left\{ y \in X \mbox{: } d^{x}(x,y)
\leq 1 \right\}$
given by
\[(\delta, \mu, x)(u,v) = \frac{1}{\mu} d(\delta^{x}_{\mu} u ,
\delta^{x}_{\mu} v) \]
Then the curve $\displaystyle \mu> 0 \mapsto \mathbb{P}^{x}(\mu) =
 [\bar{B}_{d^{x}}(x,1), (\delta, \mu, x), x]$ admits an extension by
 continuity to a   metric profile, by setting $\displaystyle
 \mathbb{P}^{x}(0) =  [\bar{B}_{d^{x}}(x,1), d^{x}, x]$.
  More precisely we have the following
 estimate:
 \begin{gather*}
 d_{GH}\left(  [\bar{B}_{d^{x}}(x,1), (\delta, \varepsilon \mu, x), x] ,
 \left[
 \bar{B}_{\frac{1}{\varepsilon}(\delta^{x}, \mu, x)}(x,1),
 \frac{1}{\varepsilon}(\delta^{x}, \mu, x), x \right] \right) \, = \\
  = \,
 \mathcal{O}(\varepsilon \mu) +
 \frac{1}{\varepsilon} \mathcal{O}(\mu) + \mathcal{O}(\mu)
 \end{gather*}
 uniformly with respect to $x$ in compact set.
 \end{thm}

\subsection{Tangent bundle of a dilatation structure}
\label{induced}

The following two theorems describe the most important metric and algebraic
properties of a dilatation structure. As presented
here these are condensed statements, available in full length as theorems 7, 8,
10 in \cite{buligadil1}.

\begin{thm}
Let $(X,d,\delta)$ be a  dilatation structure. Then the metric space $(X,d)$
admits a metric tangent space at $x$, for any point $x\in X$.
More precisely we have  the following limit:
$$\lim_{\varepsilon \rightarrow 0} \ \frac{1}{\varepsilon} \sup \left\{  \mid d(u,v) - d^{x}(u,v) \mid \mbox{ : } d(x,u) \leq \varepsilon \ , \ d(x,v) \leq \varepsilon \right\} \ = \ 0 \ .$$
\label{thcone}
\end{thm}

\begin{thm}
Let $(X,d,\delta)$ be a dilatation structure. Then for any $x \in X$ the triple
 $\displaystyle (U(x), \Sigma^{x}, \delta^{x}, d^{x})$ is a normed local
 conical group. This means:
 \begin{enumerate}
 \item[(a)]  $\displaystyle \Sigma^{x}$ is a local group operation on $U(x)$,
 with $x$ as neutral element and $\displaystyle \, inv^{x}$ as the inverse element
 function;
  \item[(b)] the distance $\displaystyle d^{x}$ is left invariant with respect
  to the group operation from point (a);
 \item[(c)] For any $\varepsilon \in \Gamma$, $\nu(\varepsilon) \leq 1$, the
 dilatation $\displaystyle \delta^{x}_{\varepsilon}$ is an automorphism with
 respect to the group operation from point (a);
 \item[(d)] the distance $d^{x}$ has the cone property with
respect to dilatations: foar any $u,v \in X$ such that $\displaystyle d(x,u)\leq 1$ and
$\displaystyle d(x,v) \leq 1$  and all $\mu \in (0,A)$ we have:
$$d^{x}(u,v) \ = \ \frac{1}{\mu} d^{x}(\delta_{\mu}^{x} u , \delta^{x}_{\mu} v)
 \quad .$$
 \end{enumerate}
\label{tgene}
\end{thm}

The conical group $\displaystyle (U(x), \Sigma^{x}, \delta^{x})$ can be regarded as the tangent space
of $(X,d, \delta)$ at $x$. Further will be denoted by:
$\displaystyle T_{x} X =  (U(x), \Sigma^{x}, \delta^{x})$.

The following is corollary 4.7 \cite{buligacont}.

\begin{thm}
Let $(X,d,\delta)$ be a dilatation structure.
Then for any $x \in X$ the local group
 $\displaystyle (U(x), \Sigma^{x})$ is locally a simply connected Lie group
 whose Lie algebra admits a positive graduation (a Carnot group).
 \label{cortang}
\end{thm}

  \section{Examples of dilatation structures}

In this section we give some examples of dilatation structures, which share some common features. There are other examples, typically coming from iterated functions systems, which will be presented
in another paper.

The first example is known to everybody: take $\displaystyle (X,d)  =  ( \mathbb{R}^{n}, d_{E})$, with usual (euclidean) dilatations $\displaystyle \delta^{x}_{\varepsilon}$, with:
$$d_{E}(x,y) = \| x-y \| \ \ \ , \ \ \delta_{\varepsilon}^{x} y \ = \ x + \varepsilon (y-x) \ .$$
Dilatations are defined everywhere. There are few things to check:  axioms 0,1,2 are obviously true. For axiom A3, remark that
for any $\varepsilon > 0$, $x,u,v \in X$ we have:
$$\frac{1}{\varepsilon} d_{E}(\delta^{x}_{\varepsilon} u , \delta^{x}_{\varepsilon} v ) \ = \ d_{E}(u,v) \ , $$
therefore for any $x \in X$ we have $\displaystyle d^{x} = d_{E}$.

Finally, let us check the axiom A4. For any $\varepsilon > 0$ and $x,u,v \in X$ we have
$$\delta_{\varepsilon^{-1}}^{\delta_{\varepsilon}^{x} u} \delta_{\varepsilon}^{x} v \ = \
x + \varepsilon  (u-x) + \frac{1}{\varepsilon} \left( x+ \varepsilon(v-x) - x - \varepsilon(u-x) \right) \ = \ $$
$$ = \ x + \varepsilon  (u-x) + v - u$$
therefore this quantity converges to
$$x + v - u \ = \ x + (v - x) - (u - x)$$
as $\varepsilon \rightarrow 0$. The axiom A4 is verified.

We continue further with less obvious examples.

\subsection{Riemannian manifolds}

Take now $\phi: \mathbb{R}^{n} \rightarrow \mathbb{R}^{n}$ a bi-Lipschitz diffeomorphism. Then
we can define the dilatation structure: $X \ = \ \mathbb{R}^{n}$,
$$d_{\phi}(x,y) \ = \ \| \phi(x) - \phi(y) \| \ \ \ , \ \ \delta_{\varepsilon}^{x} y \ = \ x+ \varepsilon(y-x) \ ,$$
or the equivalent dilatation structure:
$X \ = \ \mathbb{R}^{n}$,
$$d_{\phi}(x,y) \ = \ \| x - y \| \ \ \ , \ \ \delta_{\varepsilon}^{x} y \ = \ \phi^{-1} \left( \phi(x) + \varepsilon (\phi(y) - \phi(x))\right) .$$
In this example (look at its first  version) the distance $\displaystyle d_{\phi}$ is not equal to $\displaystyle d^{x}$. Indeed, a direct calculation shows that
$$d^{x}(u,v) \ = \ \| D\phi(x) (v-u) \| \ .$$
The axiom A4 gives the same result as previously.

Because dilatation structures are defined by local requirements, we can easily define dilatation
structures on riemannian manifolds, using particular atlases of the manifold and the riemannian
distance (infimum of length of curves joining two points).  This class of examples covers all
dilatation structures used in
differential geometry. The axiom A4 gives an operation of addition of vectors
in the tangent space  (compare with Bella\"{\i}che \cite{bell} last section).

\subsection{Snowflakes}

The next example is a snowflake variation of the euclidean case: $X =  \mathbb{R}^{n}$ and  for any $a \in (0,1]$ take
$$ d_{a}(x,y) \ = \ \| x-y \|^{\alpha} \ \ \ ,  \ \ \delta^{x}_{\varepsilon} y \ =  \ x + \varepsilon^{\frac{1}{a}} (y - x) \ .$$
We leave to the reader to verify the axioms.

More general, if $(X,d,\delta)$ is a dilatation structure then $(X,d_{a}, \delta(a))$ is also a dilatation structure, for any $a \in (0,1]$, where
$$d_{a}(x,y) \ = \ \left( d(x,y) \right)^{a} \ \ , \ \ \delta(a)_{\varepsilon}^{x}\ = \ \delta^{x}_{\varepsilon^{\frac{1}{a}}} \ .$$

\subsection{Nonstandard dilatations in the euclidean space}

Take $\displaystyle X = \mathbb{R}^{2}$ with the euclidean distance. For any $z \in \mathbb{C}$ of the
form $z= 1+ i \theta$ we define dilatations
$$\delta_{\varepsilon} x = \varepsilon^{z} x  \ .$$
It is easy to check that $(X,\delta, +, d)$ is a conical group, equivalenty that the dilatations
$$\delta^{x}_{\varepsilon} y = x + \delta_{\varepsilon} (y-x)  \ .$$
form a dilatation structure with the euclidean distance.

Two such dilatation structures (constructed with the help of complex numbers
$1+ i \theta$ and $1+ i \theta'$) are equivalent if and only if $\theta = \theta'$.

There are two other surprising properties of these dilatation structures. The first is that if $\theta \not = 0$ then there are no non trivial Lipschitz curves in $X$ which are differentiable almost everywhere.
The second property is that any holomorphic and Lipschitz function from $X$ to $X$ (holomorphic in the
usual sense on $X = \mathbb{R}^{2} = \mathbb{C}$) is differentiable almost everywhere, but there are
Lipschitz functions from $X$ to $X$ which are not differentiable almost everywhere (suffices to take a
$\displaystyle \mathcal{C}^{\infty}$ function from  $\displaystyle \mathbb{R}^{2}$ to $\displaystyle \mathbb{R}^{2}$ which is not holomorphic).

\section{Sub-riemannian dilatation structures}

To any normal frame of a regular sub-riemannian manifold we associate a
dilatation structure. (Technically this is a dilatation structure defined only
locally, as in the case of riemannian manifolds.)

\begin{defi}
To any  normal  frame
$\displaystyle \left\{ X_{1}, ..., X_{n} \right\}$ of a regular sub-riemannian
manifold $(M,D,g)$ we associate the  dilatation structure $(M,d, \delta)$
 defined by: $d$ is  the
Carnot-Carath\'eodory distance, and for any point $x \in M$
and any $\varepsilon \in (0,+\infty)$ (sufficiently small if necessary),
the dilatation $\displaystyle \delta^{x}_{\varepsilon}$ is given by:
$$\delta^{x}_{\varepsilon} \left(\exp\left( \sum_{i=1}^{n} a_{i} X_{i} \right)(x)\right) \  = \
\exp\left( \sum_{i=1}^{n} a_{i} \varepsilon^{deg X_{i}}  X_{i} \right)(x)$$
\label{dsrdil}
\end{defi}

We shall prove that $(M,d, \delta)$ is indeed a dilatation structure. This
allows us to get the main results concerning the infinitesimal geometry of a
regular sub-riemannian manifold, as particular cases of theorems \ref{thcone1},
\ref{thcone}, \ref{tgene} and \ref{cortang}.

We only have to prove axioms A3 and A4 of dilatation structures. We do this in the
next two theorems. Before this let us decribe what we mean by "sufficiently closed".

\begin{convention}
Further we shall say that  a property
$\displaystyle \mathcal{P}(x_{1}$, $\displaystyle x_{2}$,
$\displaystyle x_{3}, ...)$ holds for
$\displaystyle x_{1}, x_{2}, x_{3},...$ sufficiently closed if for any
compact, non empty set $K \subset X$, there is a positive constant $C(K)> 0$
such that $\displaystyle \mathcal{P}(x_{1},x_{2},x_{3}, ...)$ is true for any
$\displaystyle x_{1},x_{2},x_{3}, ... \in K$ with
$\displaystyle d(x_{i}, x_{j}) \leq C(K)$.
\end{convention}

In the following we prove a result similar to Gromov local approximation theorem \cite{gromovsr}, p. 135, 
or to Bella\"{\i}che theorem 7.32 \cite{bell}. Note however that here we take as a hypothesis  the existence 
of a normal frame. 

\begin{thm}
Consider $X_{1}, ... , X_{n}$ a normal frame and the associated dilatations provided by definition
\ref{dsrdil}. Then axiom A3 of dilatation structures is satisfied, that is the limit
$$\lim_{\varepsilon \rightarrow 0} \frac{1}{\varepsilon} \, d\left(
\delta^{x}_{\varepsilon} u, \delta^{x}_{\varepsilon} v \right) \ = \
d^{x}(u,v)$$
exists and it uniform with respect to x,u,v sufficiently closed.
\label{thma3}
\end{thm}

\paragraph{Proof.}
Let $x, u, v \in M$ be sufficiently closed. We write
$$u \ = \ \exp \left( \sum_{1}^{n} u_{i} X_{i} \right)(x) \quad \quad , \quad
v \ = \ \exp \left( \sum_{1}^{n} v_{i} X_{i} \right)(x)$$
we compute, using definition \ref{dsrdil}:
$$\frac{1}{\varepsilon} \, d\left(
\delta^{x}_{\varepsilon} u, \delta^{x}_{\varepsilon} v \right) \ = \
\frac{1}{\varepsilon} \, d\left(
\delta^{x}_{\varepsilon} \exp \left( \sum_{1}^{n} u_{i} X_{i} \right)(x),
\delta^{x}_{\varepsilon} \exp \left( \sum_{1}^{n} v_{i} X_{i} \right)(x) \right)
\ =$$
$$= \ \frac{1}{\varepsilon} \, d\left(
 \exp \left( \sum_{1}^{n} \varepsilon^{deg \, X_{i}} u_{i} X_{i} \right)(x),
\exp \left( \sum_{1}^{n} \varepsilon^{deg \, X_{i}}  v_{i} X_{i} \right)(x) \right)
\ = \ A_{\varepsilon}$$
Let us denote by $\displaystyle u_{\varepsilon} = \exp \left( \sum_{1}^{n}
\varepsilon^{deg \, X_{i}} u_{i} X_{i} \right)(x)$. Use the first part of the
property (b), definition \ref{defnormal} of a normal system, to write further:
$$A_{\varepsilon} \ = \ \frac{1}{\varepsilon} \, d \left( u_{\varepsilon},
\exp \left( \sum_{1}^{n} P_{i}(\varepsilon^{deg \, X_{j}} v_{j}, \varepsilon^{deg \, X_{k}} u_{k}, x)
X_{i}\right) (u_{\varepsilon}) \right) \ = $$
$$= \ \frac{1}{\varepsilon} \, d \left( u_{\varepsilon},
\exp \left( \sum_{1}^{n} \varepsilon^{deg \, X_{i}} \left( \varepsilon^{-deg \,
X_{i}} \, P_{i}(\varepsilon^{deg \, X_{j}} v_{j}, \varepsilon^{deg \, X_{k}} u_{k}, x) \right)
X_{i}\right) (u_{\varepsilon}) \right) $$
We make a final notation: for any $i=1, ..., n$
$$a_{i}^{\varepsilon} \ = \ \varepsilon^{-deg \,
X_{i}} \, P_{i}(\varepsilon^{deg \, X_{j}} v_{j}, \varepsilon^{deg \, X_{k}} u_{k}, x)$$
thus we have:
$$\frac{1}{\varepsilon} \, d\left(
\delta^{x}_{\varepsilon} u, \delta^{x}_{\varepsilon} v \right) \ = \
\frac{1}{\varepsilon} \, d \left( u_{\varepsilon},
\exp \left( \sum_{1}^{n} \varepsilon^{deg \, X_{i}} a^{\varepsilon}_{i}
X_{i}\right) (u_{\varepsilon}) \right)$$
By the second part of property (b), definition \ref{defnormal}, the vector
$\displaystyle a^{\varepsilon} \in \mathbb{R}^{n}$ converges to a finite value $\displaystyle a^{0} \in 
\mathbb{R}^{n}$, as $\varepsilon \rightarrow 0$,  uniformly with respect
to $x, u, v$ in compact set. In the same time $\displaystyle u_{\varepsilon}$ converges to $x$, as 
$\varepsilon \rightarrow 0$.  The proof ends by using
property (a), definition \ref{defnormal}. Indeed, we shall use the key assumption of uniform convergence. 
With the notations from definition \ref{defnormal}, for fixed $\eta > 0$ the term 
$$B(\eta, \varepsilon) = \frac{1}{\varepsilon} \, d \left( u_{\eta}, \exp \left( \sum_{1}^{n} \varepsilon^{deg \, X_{i}} a^{\eta}_{i}
      X_{i}\right) (u_{\eta}) \right)$$
converges to a real number $\displaystyle    A(u_{\eta}, a_{\eta})$ as  $\varepsilon \rightarrow 0$, uniformly with respect to 
$\displaystyle u_{\eta}$ and $\displaystyle a_{\eta}$. Since  $\displaystyle u_{\eta}$    converges to $x$  and 
 $\displaystyle a_{\eta}$ converges to   $\displaystyle a^{0}$ as $\eta \rightarrow 0$, by the uniform convergence assumption 
in (a), definition \ref{defnormal} we get that 
$$\lim_{\varepsilon \rightarrow 0}  \frac{1}{\varepsilon} \, d\left(  \delta^{x}_{\varepsilon} u, \delta^{x}_{\varepsilon} v \right) \  = \ 
\lim_{\eta \rightarrow 0} A(u_{\eta}, a_{\eta})  \ = \   A(x, a^{0})$$
The proof is done.                                 $\quad \square$

In the next Theorem we prove that axiom A4 of dilatation structures is satisfied. The referee 
informed us that Theorem \ref{pcomp3} also follows from results of Vodopyanov 
and Karmanova \cite{vodokar}, quoted in \cite{vodopis2} p. 267; a  complete version of
this result will apear in a work by Karmanova and Vodopyanov ``Geometry of Carnot-Carath\'eodory spaces, 
differentiability and coarea formula'' in the book ``Analysis and Mathematical Physics'', Birchh\"{a}user 2008.

\begin{thm}
Consider $X_{1}, ... , X_{n}$ a normal frame and the associated dilatations provided by definition
\ref{dsrdil}. Then axiom A4 of dilatation structures is satisfied:
 as $\varepsilon$ tends to $0$ the quantity
$$\Delta^{x}_{\varepsilon}(u,v) \ = \ \delta_{\varepsilon^{-1}}^{\delta_{\varepsilon}^{x}
u} \circ \delta_{\varepsilon}^{x}(v)$$
converges, uniformly  with respect to $x,u,v$ sufficiently closed. \label{pcomp3}
\end{thm}

\paragraph{Proof.}
We shall use the notations from definition \ref{defadapt}, \ref{defnormal},
\ref{dsrdil}.

Let $x, u, v \in M$ be sufficiently closed. We write
$$u \ = \ \exp \left( \sum_{1}^{n} u_{i} X_{i} \right)(x) \quad \quad , \quad
v \ = \ \exp \left( \sum_{1}^{n} v_{i} X_{i} \right)(x)$$
We compute now $\displaystyle \Delta^{x}_{\varepsilon}(u,v)$:
$$\Delta^{x}_{\varepsilon}(u,v) \ = \ \delta_{\varepsilon^{-1}}^{\exp \left(
\sum_{1}^{n} \varepsilon^{deg \, X_{i}} u_{i} X_{i} \right)(x)} \exp \left(
\sum_{1}^{n} \varepsilon^{deg \, X_{i}} v_{i} X_{i} \right)(x)$$
Let us denote by $\displaystyle u_{\varepsilon} = \delta^{x}_{\varepsilon} u$. Thus
we have
$$\Delta^{x}_{\varepsilon}(u,v) \ = \ \delta_{\varepsilon^{-1}}^{ u_{\varepsilon}}
\exp \left(
\sum_{1}^{n} \varepsilon^{deg \, X_{i}} v_{i} X_{i} \right)(x)$$
We use the first part of the property (b), definition \ref{defnormal}, in order to
write
$$\exp \left(
\sum_{1}^{n} \varepsilon^{deg \, X_{i}} v_{i} X_{i} \right)(x) \ = \
\exp \left( \sum_{1}^{n} P_{i}(\varepsilon^{deg \, X_{j}} v_{j}, \varepsilon^{deg \, X_{k}} u_{k}, x)
X_{i}\right) (u_{\varepsilon}) $$
We finish the computation:
$$\Delta^{x}_{\varepsilon}(u,v) \ = \ \exp \left( \sum_{1}^{n} \varepsilon^{- \, deg \, X_{i}} \, P_{i}(\varepsilon^{deg \, X_{j}} v_{j}, \varepsilon^{deg \, X_{k}} u_{k}, x)
X_{i}\right) (u_{\varepsilon})$$
As $\varepsilon$ goes to $0$ the point  $\displaystyle u_{\varepsilon}$ converges to
$x$ uniformly with respect to $x,u$ sufficiently closed (as a corollary of the
previous theorem, for example). The proof therefore ends by invoking the second
part of the property (b), definition \ref{defnormal}. $\quad \square$

\end{document}